\documentclass[12pt]{article}
\usepackage{amsmath,amsfonts,latexsym,amsthm,amssymb,mathabx}
\newcommand{\Rn}{\mathbb R^n}

\DeclareMathOperator{\R}{Re}
\DeclareMathOperator{\I}{Im}
\newcommand{\DK}{\mathbb D_{(k)}}

\newcommand{\K}{\mathcal K}

\newcommand{\CBF}{\mathcal C\mathcal B\mathcal F}
\begin{document}
\newtheorem*{prop}{Proposition}
\newtheorem{lem}{Lemma}[section]
\newtheorem*{teo}{Theorem}
\pagestyle{plain}
\title{ Growth Equation of the General Fractional Calculus}

\author{ \textbf{Anatoly N.
Kochubei}\\
Institute of Mathematics,\\
National Academy of Sciences of Ukraine, \\
Tereshchenkivska 3, \\
Kyiv, 01024 Ukraine\\
Email: kochubei@imath.kiev.ua
\and
\textbf{Yuri Kondratiev}\\
Department of Mathematics, University of Bielefeld, \\
D-33615 Bielefeld, Germany,\\
Email: kondrat@math.uni-bielefeld.de}

\date{}
\maketitle

\begin{abstract}

We consider the Cauchy problem $(\DK u)(t)=\lambda u(t)$, $u(0)=1$, where $\DK$ is the
general convolutional  derivative introduced in the paper (A. N. Kochubei, Integral Equations Oper. Theory {\bf 71} (2011), 583--600), $\lambda >0$. The solution is a generalization of the function $t\mapsto E_\alpha (\lambda t^\alpha)$ where $0<\alpha <1$, $E_\alpha$ is the Mittag-Leffler function. The asymptotics of this solution, as $t\to \infty$, is studied.

\medskip{}

\end{abstract}

\section{Introduction}

In several models of dynamics  of complex systems the time evolution for observed quantities  has exponential asymptotics
of two possible types. In the simplest cases these asymptotics are related with the solutions to  the equations
$$
u'(t)=zu(t),\quad t>0;\quad u(0)=1,
$$
where we will consider positive and negative $z$ separately.
For $z<0$ (the relaxation equation) the solution   decays  to zero. In particular models as, e.g., Glauber
stochastic dynamics in the continuum, this corresponds to an exponential  convergence to an equilibrium,
see \cite{KKuM}.  The case $z>0$  also may appear in applications. We can mention
the contact model in the continuum where  for the mortality  below a critical value the density
of population will grow exponentially fast \cite{KS,KKuP}, as well as models of economic growth.

On the other hand, the observed behavior of specific physical and biological systems show an emergence
of other time asymptotics which may be far from the exponential  decay or growth.  An attempt to obtain
other relaxation characteristics is related with a use of generalized time derivatives in dynamical equations (see \cite{KL1,KL2}).
In this way we may produce a wide spectrum of possible asymptotics to reflect a demand coming from
applications \cite{TM}.

The general fractional calculus introduced in \cite{K2011} is based on a version of the fractional derivative, the differential-convolution operator
$$
(\DK u)(t)=\frac{d}{dt}\int\limits_0^tk(t-\tau )u(\tau )\,d\tau -k(t)u(0)
$$
where $k$ is a nonnegative locally integrable function satisfying additional assumptions, under which

(A) The Cauchy problem
\begin{equation}
\label{1}
(\DK u)(t)=-\lambda u(t),\quad t>0;\quad u(0)=1,
\end{equation}
where $\lambda >0$, has a unique solution which is completely monotone;

(B) The Cauchy problem
\begin{equation*}
(\DK w)(t,x)=\Delta w(t,x),\quad t>0,\ x\in \Rn;\quad w(0,x)=w_0(x),
\end{equation*}
is solvable (under appropriate conditions for $w_0$) and possesses a fundamental solution, a kernel with the property of a probability density.

A class of functions $k$, for which (A) and (B) hold, was found in \cite{K2011} and is described below. The simplest example is
\begin{equation}
\label{2}
k(t)=\frac{t^{-\alpha }}{\Gamma (1-\alpha )},\quad t>0,
\end{equation}
where $0<\alpha <1$, and for this case $\DK$ is the Caputo-Djrbashian fractional derivative $\mathbb D^{(\alpha )}$. Another subclass is the one of distributed order derivatives; see \cite{K2008} for the details.

Note that for the case where $k$ has the form (\ref{2}), the solution of (\ref{1}) is $u(t)=E_\alpha (-\lambda t^\alpha)$ where $E_\alpha$ is the Mittag-Leffler function; see \cite{KST}, Lemma 2.23 (page 98). This solution has a slow decay at infinity, due to the asymptotic property of the Mittag-Leffler function; see \cite{GKMR}.  Note that using particular classes of fractional derivatives we observe
several specific asymptotics for the solution of the equation (\ref{1}) with $\lambda >0$. Some results in this direction were already obtained
in \cite{K2008,K2009,KK}. A more detailed analysis of this problem will be performed in a forthcoming paper.

In this paper, we consider the Cauchy problem with the opposite sign in the right-hand side, that is
\begin{equation}
\label{3}
(\DK u)(t)=\lambda u(t),\quad t>0;\quad u(0)=1;
\end{equation}
as before, $\lambda >0$. In the case (\ref{2}), we have $u(t)=E_\alpha (\lambda t^\alpha )$ (see \cite{KST}, Lemma 2.23 (page 98)), and, due to the well-known asymptotics of $E_\alpha$ \cite{GKMR}, this is a function of exponential growth. The existence and uniqueness of an absolutely continuous solution of (\ref{3}) follows from the results of \cite{Sin} dealing with more general nonlinear equations. Here we study the asymptotic behavior of the solution of (\ref{3}). Functions of this kind can be useful for fractional macroeconomic models with long dynamic memory; see \cite{TT} and references therein. Let us explain this in a little greater detail.

In modern macroeconomics, the most important are so-called growth models,
which in the mathematical sense are reduced (for linear models)
to the equation $u'(t)=\lambda u(t) +f(t)$ with $\lambda >0$.
In economics, an important role is played by processes with a distributed lag
starting with Phillips' works \cite{Ph} (see also \cite{Al}) and long memory, starting with Granger's work \cite{Gr} (see also \cite{Bai}).

If we assume the presence of effects of distributed lag (time delay) or fading memory
in economic processes, then the fractional generalization of the linear classical growth models
can be described by the fractional differential equation
$D^{\alpha} u(t) = \lambda u(t) +f(t)$ with $\lambda >0, \alpha>0$.
The fractional generalizations of well-known economics models
have been first proposed for the Caputo-Djrbashian fractional derivative $D^\alpha$.
Solving the problem in a more general case will allow us to
describe accurately the conditions on the operator kernels (the memory functions),
under which equations for models of economic growth with memory have solutions.

In general fractional calculus, which was proposed in \cite{K2011} (see also \cite{KL1}),
the case $\lambda >0$ is not considered.
The growth equation was considered in \cite{K2008} for the special case of a
distributed order derivative, where it was proved that
a smooth solution exists and is monotone increasing.

In this article we propose correct mathematical statements
for growth models with memory in more general case,
for the general fractional derivative $\DK$
with respect to the time variable. Their application
can be useful for mathematical economics for the description of processes with long memory and distributed lag.

Note that the technique used below was developed initially in \cite{KK} for the use in the study of intermittency in fractional models of statistical mechanics.

\section{Preliminaries}

Our conditions regarding the function $k$ will be formulated in terms of its Laplace transform
\begin{equation}
\label{4}
\K (p)=\int\limits_0^\infty e^{-pt}k(t)\,dt.
\end{equation}
Denote $\Phi (p)=p\K (p)$.

We make the following assumptions leading to (A) and (B) (see \cite{K2011}).

\begin{description}
\item[(*)] The Laplace transform (\ref{4}) exists for all positive numbers $p$. The function $\K$ belongs to the Stieltjes class $\mathcal S$, and
    \begin{equation}
    \label{5}
    \K (p)\to \infty ,\text{ as $p\to 0$};\quad \K (p)\to 0,\text{ as $p\to \infty$};
    \end{equation}
     \begin{equation}
     \label{6}
    p\K (p)\to 0,\text{ as $p\to 0$};\quad p\K (p)\to \infty,\text{ as $p\to \infty$};
    \end{equation}
\end{description}

Recall that the Stieltjes class consists of the functions $\psi$ admitting the integral representation
$$
\psi (z)=\frac{a}{z}+b +\int\limits_0^\infty \frac1{z+t}\sigma (dt)
$$
where $a,b\ge 0$, $\sigma$ is a Borel measure on $[0,\infty )$, such that
\begin{equation}
\label{7}
\int\limits_0^\infty (1+t)^{-1}\sigma (dt)<\infty .
\end{equation}
For a detailed exposition of the theory of Stieltjes functions including properties of the measure $\sigma$ see \cite{SSV}, especially Chapters 2 and 6.

In particular, for the Stieltjes function $\K$, the limit conditions (\ref{5}), (\ref{6}) imply the representation
\begin{equation}
\label{8}
\K (p)=\int\limits_0^\infty \frac1{z+t}\sigma (dt).
\end{equation}
We can also write \cite{K2011} that
$$
k(s)=\int\limits_0^\infty e^{-ts}\,\sigma (dt),\quad 0<s<\infty.
$$

The function $\Phi$ belongs to the class $\CBF$ of complete Bernstein functions, a subclass of the class $\mathcal B\mathcal F$ of Bernstein functions. Recall that a function $f:\  (0,\infty )\to \mathbb R$ is called a Bernstein function, if $f\in C^\infty$, $f(z)\ge 0$ for all $z>0$, and
$$
(-1)^{n-1}f^{(n)}(z)\ge 0\quad \text{for all $n\ge 1,z>0$}.
$$
so that the derivative of $f$ is completely monotone. A function $f$ belongs to $\CBF$, if it has an analytic continuation to the cut complex plane $\mathbb C\setminus (-\infty ,0]$, such that $\I z\cdot \I f(z)\ge 0$, and there exists the real limit
$$
f(0+)=\lim\limits_{(0,\infty )\ni z\to 0}f(z).
$$

Both the classes $\mathcal B\mathcal F$ and $\CBF$ admit equivalent descriptions in terms of integral representations; see \cite{SSV}.

Below we will need the following inequality for complete Bernstein functions (Proposition 2.4 in \cite{BCT}), valid, in particular, for the function $\Phi$. For any $p$ outside the negative real semi-axis, we have
\begin{equation}
\label{9}
\sqrt{\frac{1+\cos \varphi}2}\Phi (|p|)\le |\Phi (p)|\le \sqrt{\frac2{1+\cos \varphi}}\Phi (|p|),\quad \varphi =\arg p.
\end{equation}

Solutions of the Cauchy problem (\ref{3}) and a similar problem with the classical first order derivative are connected by the subordination identity (see \cite{K2011}; for the case of the Caputo-Djrbashian derivative see \cite{Ba}), an integral transformation with the kernel $G(s,t)$ constructed as follows.

Consider the function
\begin{equation}
\label{10}
g(s,p)=\K (p)e^{-s\Phi (p)},\quad s>0,p>0.
\end{equation}
It is proved \cite{K2011} that $g$ is a Laplace transform in the variable $t$ of the required kernel $G(s,t)$, that is
$$
g(s,p)=\int\limits_0^\infty e^{-pt}G(s,t)\,dt.
$$
$G$ is nonnegative, and
$$
\int\limits_0^\infty G(s,t)\,ds=1 \quad \text{for each\ $t$}.
$$

\section{Cauchy problem for the growth equation}

Let us consider the Cauchy problem (\ref{3}). If $u_\lambda (t)$ is its solution whose Laplace transform $\widetilde{u_\lambda}(p)$ exists for some $p$, then it follows from properties of the Laplace transform \cite{Do} that
$$
\Phi (p)\widetilde{u_\lambda}(p)-\lambda \widetilde{u_\lambda}(p)=\K (p),
$$
hence
\begin{equation}
\label{11}
\widetilde{u_\lambda}(p)=\frac{\K (p)}{\Phi (p)-\lambda},\quad \text{if $\Phi (p)>\lambda$.}
\end{equation}

On the other hand, consider the function
\begin{equation}
\label{12}
E(t,\lambda )=\int\limits_0^\infty e^{\lambda s}G(s,t)\,ds,\quad t>0.
\end{equation}
The existence of the integral in (12) for almost all $t>0$ is, by the Fubini-Tonelli theorem, a consequence of the absolute convergence of the repeated integral
$$
\int\limits_0^\infty e^{\lambda s}\,ds\int\limits_0^\infty e^{-pt}G(s,t)\,dt=\int\limits_0^\infty e^{\lambda s}g(s,\lambda )\,ds=\frac{\K (p)}{\Phi (p)-\lambda }
$$
where $p>0$ is such that $\Phi (p)>\lambda $.

The above calculation shows that $E(t,\lambda )=u_\lambda (t)$, the solution of (\ref{3}), and the identity (\ref{12}) provides an integral representation of this solution.

A more detailed analysis of its properties is based on the analytic properties of the Stieltjes function $\K$ or, equivalently, of the complete Bernstein function $\Phi$; in particular, we use the representation
\begin{equation}
\label{13}
\Phi (p)=\int\limits_0^\infty \frac{p}{p+t}\,\sigma (dt),
\end{equation}
which follows from (\ref{7}). The measure $\sigma$ satisfies (\ref{8}).

Since $\Phi$ is a Bernstein function, its derivative $\Phi'$ is completely monotone. By our assumptions, $\Phi$ is not a constant function, so that $\Phi'$ is not the identical zero. It follows from Bernstein's description of completely monotone functions that $\Phi'(p)\ne 0$ for any $p>0$ (see Remark 1.5 in \cite{SSV}). Therefore $\Phi$ is strictly monotone, and for each $z>0$, there exists a unique $p_0=p_0(z)>0$, such that $\Phi (p_0)=z$. The inequality $\Phi (p)>z$ is equivalent to the inequality $p>p_0(z)$. Since $\Phi$, as a complete Bernstein function, preserves the open upper and lower half-planes (in fact, this follows from (\ref{13}), we have $\Phi (p)\ne z$ for any nonreal $p$.

It is proved in \cite{KK} that the function $p_0(z)$, $z>0$, is strictly superadditive, that is
$$
p_0(x+y)>p_0(x)+p_0(y)\quad \text{for any $x,y>0$}.
$$

\medskip
\begin{prop}
The solution $u_\lambda (t)$ of the Cauchy problem (\ref{3}) admits a holomorphic continuation in the variable $t$ to a sector $\Sigma_v=\left\{ re^{i\theta}:\ r>0,-v<\theta <v\right\}$, $0<v<\frac{\pi}2$, and
\begin{equation}
\label{14}
\sup\limits_{t\in \Sigma_v}\left| e^{-p_0t}u_\lambda (t)\right| <\infty,\quad p_0=p_0(\lambda ).
\end{equation}
\end{prop}

\medskip
{\it Proof.} It follows from (\ref{11}) and (\ref{13}) that the Laplace transform $\widetilde{u_\lambda}(p)$ is holomorphic in $p$ on any sector $p_0+\Sigma_{\rho +\frac{\pi}2}$, $0<\rho <\frac{\pi}2$. In addition,
\begin{equation}
\label{15}
\sup\limits_{p\in p_0+\Sigma_{\rho +\frac{\pi}2}}\left| (p-p_0)\widetilde{u_\lambda}(p)\right| <\infty.
\end{equation}
Now the assertion is implied by (\ref{15}) and the duality theorem for holomorphic continuations of a function and its Laplace transform; see Theorem 2.6.1 in \cite{ABHN}. $\qquad \blacksquare$

\medskip
Now we are ready to formulate and prove our main result.

\medskip
\begin{teo}
Let the assumptions $(*)$ hold, and in addition,
\begin{equation}
\label{16}
\int\limits_1^\infty \frac{ds}{s\Phi (s)}<\infty .
\end{equation}
Then
\begin{equation}
\label{17}
u_\lambda (t)=\frac{\lambda }{\Phi'(p_0(\lambda ))p_0(\lambda )}e^{p_0(\lambda )t}+o(e^{p_0(\lambda )t}),\quad t\to \infty.
\end{equation}
\end{teo}

\medskip
{\it Proof.} The representation (\ref{11}) can be written as
$$
\widetilde{u_\lambda}(p)=\frac1p\left( 1+\frac{\lambda}{\Phi (p)-\lambda}\right) .
$$
This implies the representation of $u_\lambda$ as $u_\lambda (t)=1+B_\lambda (t)$ where $B_\lambda$ has the Laplace transform
$$
\widetilde{B_\lambda}(p)=\frac{\lambda}p\cdot \frac1{\Phi (p)-\lambda},
$$
for such $p$ that $\Phi (p)>\lambda$.

Using the inequality (9) we find that $|\Phi (p)|\ge \frac1{\sqrt{2}}\Phi (|p|)$ on any vertical line $\{ p=\gamma +i\tau, \tau\in \mathbb R\}$ where $\gamma >p_0$. By our assumption (16), $\widetilde{B_\lambda}$ is absolutely integrable on such a line. In addition, it follows from (\ref{6}) that $\widetilde{B_\lambda}(p)\to 0$, as $p\to \infty$ in the half-plane $\R p>p_0$. These properties make it possible (see Theorem 28.2 in \cite{Do}) to write the inversion formula
$$
u_\lambda (t)=1+\frac{\lambda }{2\pi i}\int\limits_{\gamma -i\infty}^{\gamma +i\infty} e^{pt}\frac{dp}{p(\Phi (p)-\lambda )},\quad \gamma >p_0.
$$

Denote
$$
V(t)=1+\frac{\lambda }{2\pi i}\int\limits_{r-i\infty}^{r+i\infty} e^{pt}\frac{dp}{p(\Phi (p)-\lambda )}
$$
where $0<r<p_0$. Then
\begin{equation}
\label{18}
|V(t)| \le 1+Ce^{rt}\left| \int\limits_{-\infty}^{\infty}e^{i\tau t}\frac{d\tau}{(r+i\tau )(\Phi (r+i\tau )-\lambda )}\right| =o(e^{rt}),\quad t\to \infty,
\end{equation}
by virtue of (\ref{9}), (\ref{16}) and the Riemann-Lebesgue theorem.

On the other hand, we may write
$$
u_\lambda (t)-V(t)=\frac{\lambda }{2\pi i}\left( \int\limits_{\Lambda_+}+\int\limits_{\Lambda_0}+\int\limits_{\Lambda_-}\right) e^{pt}\frac{dp}{p(\Phi (p)-\lambda )}
$$
where the contour $\Lambda_+$ consists of the vertical rays $\{ \R p=r,\I p\ge R\}$, $\{ \R p=\gamma,\I p\ge R\}$, and the horizontal segment $\{ r\le \R p\le \gamma, \I p=R\}$ ($R>0$), $\Lambda_-$ is a mirror reflection of $\Lambda_+$ with respect to the real axis, $\Lambda_0$ is the finite rectangle consisting of the vertical segments $\{ \R p=r,|\I p|\le R\}$, $\{ \R p=\gamma ,|\I p|\le R\}$, and the horizontal segments $\{ r\le \R p\le \gamma,\I p=\pm R\}$.

We have
$$
\int\limits_{\Lambda_+}e^{pt}\frac{dp}{p(\Phi (p)-\lambda )}=0,
$$
due to the Cauchy theorem, absolute integrability of the integrand on the vertical rays (see (\ref{16})) and the estimate
$$
\left| \int\limits_{\Pi_h}e^{pt}\frac{dp}{p(\Phi (p)-\lambda )}\right| \le Ch^{-1}\to 0,\quad h\to \infty,
$$
where $\Pi_h=\{ r\le \R p\le \gamma,\I p=h\}$, $h>R$. In a similar way, we prove that
$$
\int\limits_{\Lambda_-}e^{pt}\frac{dp}{p(\Phi (p)-\lambda )}=0.
$$

Due to the inequality $\Phi'(p_0)\ne 0$, there exists a complex neighborhood $W$ of the point $\lambda =\Phi (p_0)$, in which $\Phi$ has a single-valued holomorphic inverse function $p=\psi (w)$, so that $\Phi (\psi (w))=w$ and $p_0=\psi (\lambda )$. In the above arguments, the numbers $r,\gamma ,R$ were arbitrary. Now we choose $R$ and $\gamma -r$ so small that that the curvilinear rectangle $\Phi (\Lambda_0)$ lies inside $W$. Making the change of variables $p=\psi (w)$ and using the Cauchy formula we find that
\begin{multline*}
\frac{\lambda }{2\pi i}\int\limits_{\Lambda_0}e^{pt}\frac{dp}{p(\Phi (p)-\lambda )}=\frac{\lambda }{2\pi i}\int\limits_{\Phi (\Lambda_0)}e^{\psi (w)t}\frac1{\Phi'(\psi (w))\psi (w)}\cdot \frac{dw}{w-\lambda }
=\frac{\lambda }{\Phi'(\psi (\lambda ))\psi (\lambda )}e^{\psi (\lambda )t}\\
=\frac{\lambda }{\Phi'(p_0(\lambda ))p_0(\lambda )}e^{p_0(\lambda )t}.
\end{multline*}

Together with (\ref{18}), this implies the required asymptotic relation (\ref{17}). $\qquad \blacksquare$

\medskip
{\it Examples}. 1) In the case (\ref{2}) of the Caputo-Djrbashian fractional derivative of order $0<\alpha <1$, we have $u_\lambda (t)=E_\alpha (\lambda t^\alpha )$, $\Phi (p)=p^\alpha$, and the condition (\ref{16}) is satisfied. Here the asymptotics (\ref{17}) coincides with the one given by the principal term of the asymptotic expansion of the Mittag-Leffler function. The above proof is different from the classical proof of the latter (see \cite{GKMR}).

2) Let us consider the case of a distributed order derivative with a weight function $\mu$, that is
$$
\mathbb D^{(\mu )}u(t)=\int\limits_0^1 (\mathbb D^{(\alpha )})u(t)\mu
(\alpha )\,d\alpha .
$$
Suppose that $\mu \in C^2[0,1]$, $\mu (1)\ne 0$. In this case \cite{K2008},
$$
k(s)=\int\limits_0^1\frac{s^{-\alpha }}{\Gamma (1-\alpha )}\mu
(\alpha )\,d\alpha ,\quad \Phi (p)=\int\limits_0^1 p^\alpha \mu (\alpha )\,d\alpha ,
$$
and under the above assumptions,
$$
\Phi (p)=\frac{\mu (1)p}{\log p}+O\left( p|\log p|^{-2}\right) ,\quad p\to \infty.
$$
The condition (\ref{16}) is satisfied, and our asymptotic result (\ref{17}) is applicable in this case.

\section*{Acknowledgments}
The authors are grateful to V. E. Tarasov for calling our attention to the problem studied in this paper, as well as for his advices regarding its text.

The work of the first-named author was funded in part under the budget program of Ukraine No. 6541230 ``Support to the development of priority research trends'' and under the research work "Markov evolutions in real and p-adic spaces" of the Dragomanov National Pedagogical  University of Ukraine.

\end{document}